%% revisit, march 2026, hoping to get it to arXiv : 

%% topoffile 

% This is nonpmic_pmic_2.tex: in convergence! 
% 25.09.1996: version 0.08 or so
% 29.09.1996: version 0.43 or so 
% 02.10.1996: version 0.52 or so 
% 03.10.1996: version 0.59 or so 
% 04.10.1996: version 0.67 or so 
% 07.10.1996: version 0.76 or so 
% 10.10.1996: version 0.79 or so 
% 13.10.1996: version 0.84 or so 
% 30.10.1996: version 0.90 or so? with table in place 

%% need to go from old .ps things to .pdf things to get pdftex to work
%% ps2pdf heisan.ps heisan.pdf
%% pdfcrop heisan.pdf heisan-crop.pdf 

% \input psfig.sty

\input miniltx
\input graphicx

\magnification\magstep1
% \hsize15.0truecm
% \hoffset0.75truecm 
\vsize23.5truecm 
\baselineskip14pt
\font\csc=cmcsc10
\font\bigbf=cmbx12 
\font\smallrm=cmr8 

\def\hatt{\widehat}
\def\half{\hbox{$1\over2$}}
\def\d{{\rm d}}
\def\sumin{\sum_{i=1}^n}
\def\dell{\partial} 

\def\tilda{\widetilde}

\def\E{{\rm E}}
\def\RR{\mathord{I\kern-.3em R}}
\def\Var{{\rm Var}}
\def\arr{\rightarrow} 
\def\midd{\,|\,}
\def\tr{{\rm t}}

\def\N{{\rm N}}
\def\section{\medskip}
\def\mise{\hbox{\csc mise}} 
\def\mse{\hbox{\csc mse}} 
\def\ise{\hbox{\csc ise}} 
\def\ML{maximum likelihood}

\font\smallcyr=wncyr10 at 10truept
\def\Yepa{Yepanechnikov}
\def\Yep{{K_{{\hbox{\smallcyr epan}}}}}

% FERMAT: a margin-note trick 
\def\fermat#1{\setbox0=\vtop{\hsize4.00pc
        \smallrm\raggedright\noindent\baselineskip9pt
        \rightskip=0.5pc plus 1.5pc #1}\leavevmode
        \vadjust{\dimen0=\dp0
        \kern-\ht0\hbox{\kern-4.00pc\box0}\kern-\dimen0}}

\baselineskip14pt 

\bigskip
\centerline{\bigbf Sometimes nonparametrics beat parametrics,}
\centerline{\bigbf even when the model is right} 

\medskip
\centerline{\bigbf Morten Byholt and Nils Lid Hjort} 

\medskip
\centerline{\bf University of Oslo} 

\smallskip
\centerline{\sl October 1996} 
% \centerline{\sl ... but this is only \Version, as of \today} 

\medskip
{\smallskip\narrower\noindent\baselineskip12pt 
{\csc Abstract.}
A basic issue in both teaching of and practice of statistics 
is the interplay between modelling assumptions and inference performance. 
The general message conveyed is that 
stronger assumptions lead to better 
statistical performance of the relevant estimators, 
tests and confidence intervals, provided that these assumptions hold.
On the other hand, fewer assumptions often lead to 
safer and more robust methods that are good also outside 
narrow conditions, but not quite as good 
as specialist methods that exploit such narrower conditions, 
if these are fulfilled. 

This interplay is nicely illustrated in the context of density estimation,
where parametric and nonparametric methods can be contrasted. 
The parametric ones have mean squared errors of size $O(n^{-1})$ 
in terms of sample size $n$ if the parametric model is right, 
but are not even consistent outside the model. 
The nonparametric methods are everywhere consistent 
and have mean squared errors of size $O(n^{-4/5})$
for broad classes of estimands. 

The point we are making here is that this picture 
is not universally true! We show that a simple 
kernel density estimator can perform better than 
a directly estimated parametric density on the latter's home turf,
for small sample sizes, in the sense of mean integrated squared error. 
% estimated normal density on the 
% latter method's home turf, for $n\le 15$ [xx?xx], 
Our main example is that of estimating an unknown normal density. 
% which intuitively should be an easy win for the parametric method 
In the process of developing and discussing this 
somewhat counter-intuitive and half-paradoxical example 
we touch on several tangential issues of interest, 
pertaining to exact small-sample analysis of density estimators. 

\smallskip\noindent 
{\csc Key words:} \sl 
comparison of performance; 
density estimators;
nonparametric versus parametric;
small-sample analysis 
\smallskip}

\bigskip
{\bf 1. Introduction and summary.} 
There is a plethora of examples in statistics that illustrate
the following principle: stronger assumptions about the model
lead to construction of inference methods with stronger performance,
provided the assumptions in question are in force.
The other side of the coin is equally important and 
well understood in many situations: weaker assumptions 
lead to alternative methods that perform well under a broader
range of underlying probability mechanisms than those 
that require harder assumptions, but these robust methods 
typically do not perform quite as well, under a narrower model,
as the specialist methods that exploit such assumptions. 
These principles can be formulated and exemplified 
in contexts of smaller and bigger parametric models, 
or in situations with simpler and not so simple regression structures,
and so on. They are also nicely illustrated in situations
where there are both parametric and natural nonparametric competitors
for an inference task; one expects the parametric methods to behave 
worse outside the parametric model, but one surely expects 
the parametrically based optimum methods to outperform the
nonparametric ones on the parametric model's home turf. 

The main purpose of this paper is to exhibit 
a counterexample to the belief that these principles always hold. 
We shall discuss a natural situation where a simple nonparametric 
density estimator outperforms the natural method based on normality, 
even when the normal model is perfect. 
We hasten to add that this only can happen for 
small sample sizes, however. 

Suppose that $X_1,\ldots,X_n$ are independent 
and that it is required to estimate its density $f$.
If $f$ is modelled parametrically, say $f(x)=f(x,\theta)$, 
one would typically plug in natural estimators 
for the parameters to get a density estimator of the form
$\hatt f(x)=f(x,\hatt\theta)$. In the normal case this leads to
$$\hatt f(x)=\hatt\sigma^{-1}\phi(\hatt\sigma^{-1}(x-\hatt\mu)), \eqno(1.1)$$
writing $\phi(x)$ for the standard normal density, 
where $\hatt\mu=n^{-1}\sumin X_i$ and 
$\hatt\sigma^2=(n-1)^{-1}\sumin(X_i-\hatt\mu)^2$;
variations in the choice of parameter estimates 
could of course be considered. 
If on the other hand nothing is assumed a priori
besides some smoothness, one would use a nonparametric
estimator, of which there are several well-explored classes.
The most important one is the kernel type estimator,
$$\tilda f(x)=n^{-1}\sumin K_h(X_i-x), \eqno(1.2)$$
where $K_h(u)=h^{-1}K(h^{-1}u)$ is the scaled version of 
a suitable symmetric probability density function $K$
called the kernel function $K$. 

The behaviour of such estimators has been studied extensively,
see for example recent accounts by Scott (1992) and Wand and Jones (1995). 
The primary focus has been performance for moderate to large sample sizes,
however, partly because of the availability of revealing 
large-sample approximations. In particular it is well known that 
the mean squared error (\mse) of these estimators is of size 
$O(n^{-4/5})$ whereas the parametric ones, exemplified here by (1.1),
has \mse{} rate $O(n^{-1})$. Again, our aim is to show that 
for small sample sizes the picture is reversed. 
We include two kernel functions in this little study,
namely the standard normal kernel $K=\phi$ and the 
so-called \Yepa{} kernel $\Yep(u)={3\over 2}(1-4u^2)$ on $[-\half,\half]$. 
The reason is that this kernel function is known to
minimise both the pointwise \mse{} and the integrated 
mean squared error (\mise) for large samples. 

We go through the details of parametric estimators 
in Section 2 and nonparametric estimators in Section 3 and 4. 
Then Section 5 gives some comparison results and 
offers some illustrations. Finally some supplementing comments 
and results are listed in Section~6. 

It should be emphasised that this paper is consistently 
concerned with exact finite-sample properties of both types of estimators. 
This requires more attention to technical details than for 
the more easily available and implementable methods and results 
about large-sample behaviour. We do however touch on these
easier to get results too for illustration and comparison purposes. 

\bigskip
{\bf 2. MSE and MISE for the parametric estimator.}
Here we calculate exact mean, variance and \mise{}
for the natural normality-based estimator (1.1). 

\section 
{\sl 2.1. Large-sample calculations.}
Approximations in the tradition of large-sample theorems 
are both reasonably easy to obtain and illustrative. 
% so we include them here even though our emphasis is on
% the exact finite-sample results which will follow later.
For a general parametric family $f(x,\theta)$, 
let $u(x,\theta)=\dell\log f(x,\theta)/\dell\theta$
be the score function. It is well known that 
the \ML{} estimator $\hatt\theta$ behaves asymptotically according to 
$\sqrt{n}(\hatt\theta-\theta)\arr_d \N(0,J(\theta)^{-1})$,
under model conditions, 
where $J(\theta)=\E u(X,\theta)u(X,\theta)^\tr$ is 
the Fisher information matrix. By the delta method it follows that 
$$\sqrt{n}\{f(x,\hatt\theta)-f(x,\theta)\}
	\arr_d f(x,\theta)u(x,\theta)^\tr V,
	\quad {\rm where\ }V\sim \N(0,J(\theta)^{-1}). $$
In other words, $n$ times the squared error at $x$ 
tends to a variable with variance 
$v(x,\theta)^2=f(x,\theta)^2u(x,\theta)^\tr J(\theta)^{-1}u(x,\theta)$,
and $n$ times \mise{} will under mild regularity conditions
tend to $\int v(x,\theta)^2\,\d x$. 
% \fermat{so this approximation can be compared to Morten's
% simulation figure}
For the normal case this leads to the approximation 
$$\mse(x)=\E\Bigl\{{1\over \hatt\sigma}
	\phi\Bigl({x-\hatt\mu\over \hatt\sigma}\Bigr)
	-{1\over \sigma}\phi\Bigl({x-\mu\over \sigma}\Bigr)\Bigr\}^2
	\doteq {1\over n}{1\over \sigma^2}
	\phi(y)^2\{y^2+\half(y^2-1)^2\}, \eqno(2.1)$$
in terms of $y=y_x=(x-\mu)/\sigma$, and secondly to 
$$\eqalign{\mise
&\doteq n^{-1}\sigma^{-2}\int\phi(y_x)^2\{y_x^2+\half(y_x^2-1)^2\}\,\d x \cr
&=n^{-1}\sigma^{-1}\int\phi(y)^2\{y^2+\half(y^2-1)^2\}\,\d y 
={1\over n}{1\over \sigma}{1\over 2\sqrt{\pi}}
		\Bigl({4\over 8}+{3\over 8}\Bigr). \cr} \eqno(2.2)$$
The last formula indicates that 4/7 of the total \mise{} stems from
estimation of unknown $\mu$ while estimation of unknown $\sigma$ 
is responsible for the remaining 3/7. 

\section 
{\sl 2.2. Exact MSE under normality.}
We may write $\hatt\mu=\mu+(\sigma/\sqrt{n})N$ and 
$\hatt\sigma^2=\sigma^2Z^2$, say, where $N$ is standard normal
and independent of $Z^2$, a $\chi^2_{n-1}/(n-1)$ variable. 
Thus the parametric normal density estimator can be expressed as 
$$\hatt f(x)={1\over \sqrt{2\pi}}{1\over \sigma Z}
	\exp\Bigl\{-\half{(N+\sqrt{n}y)^2\over nZ^2}\Bigr\}, 
	\quad {\rm where\ again\ }y=(x-\mu)/\sigma. $$
This makes it possible to compute at first conditional moments 
given $\hatt\sigma$. One needs the fact that 
$$\E\exp\{-\half t(N+a)^2\}=(1+t)^{-1/2}
	\exp\{-\half a^2t/(1+t)\}
	\quad {\rm for\ }t>-1, \eqno(2.3)$$
which can be verified directly (and is also obtainable 
from the moment-generating function of a non-central chi squared). Hence 
$$\E\{\hatt f(x)\midd Z\}
	={1\over \sqrt{2\pi}}{1\over \sigma}{\sqrt{n}\over (1+nZ^2)^{1/2}}
	\exp\{-\half y^2n/(1+nZ^2)\} \eqno(2.4)$$
and 
$$\E\{\hatt f(x)^2\midd Z\}
	={1\over 2\pi}{1\over \sigma^2}{\sqrt{n}\over Z(2+nZ^2)^{1/2}}
	\exp\{-y^2n/(2+nZ^2)\}. \eqno(2.5)$$
It follows that the exact mean and variance, 
and hence the exact $\mse(x)$, can be arrived at by 
numerical integration with respect to the density for $Z$, 
which is of the form $g_n(z)=\gamma_{n-1}((n-1)z^2)2(n-1)z$, 
where $\gamma_{n-1}(\cdot)$ is the density of a $\chi^2_{n-1}$. 
Note that $Z$ for reasonably sized $n$ is quite close to 1 
with high probability. 

\section
{\sl 2.3. Exact MISE under normality.}
Turning next to the exact \mise, note first the technical fact that   
$$\int{1\over \sigma_1}\phi\Bigl({x-\mu_1\over \sigma_1}\Bigr)
	 {1\over \sigma_2}\phi\Bigl({x-\mu_2\over \sigma_2}\Bigr)\,\d x
        ={1\over (\sigma_1^2+\sigma_2^2)^{1/2}}
	 \phi\Bigl({\mu_1-\mu_2\over (\sigma_1^2+\sigma_2^2)^{1/2}}\Bigr); $$
a simple demonstration involves noting that $\int f_1f_2\,\d x$ 
is the density at zero for a variable of the form $Y_1-Y_2$,
where $Y_1$ and $Y_2$ are independent from respectively $f_1$ and $f_2$.
This leads to 
$$\int(\hatt f-f)^2\,\d x
=\phi(0)\Bigl[{1\over \sqrt{2}\hatt\sigma}+{1\over \sqrt{2}\sigma}
	-2{1\over (\hatt\sigma^2+\sigma^2)^{1/2}}
	 \exp\Bigl\{-\half{(\hatt\mu-\mu)^2\over 
	 \hatt\sigma^2+\sigma^2}\Bigr\}\Bigr]. $$
Now use $\E\exp(-\half t\chi^2_1)=(1+t)^{-1/2}$ 
for the moment-generating function of a squared standard normal,
a special case of (2.3). One then finds 
$$\mise=n^{-1}\sigma^{-1}(2\sqrt{\pi})^{-1}a_n, \eqno(2.6)$$
where 
$$a_n=n\,\E\Bigl[1+{1\over Z}-2\Bigl\{{2n\over 1+n(1+Z^2)}\Bigr\}^{1/2}
	\Bigr]. $$
The result is given in this form since $a_n$ is stable as $n$ grows.
Note that the formula (2.6) is exact, 
with a value for $a_n$ that can be computed via 
numerical integration using the density $g_n(z)$ for $Z$ again
(the $1/Z$ term has exact mean equal to 
$\{\half(n-1)\}^{1/2}\Gamma(\half(n-2))/\Gamma(\half(n-1))$, 
but the second term does not have a closed-form expectation formula). 
The $a_n$ sequence goes to $7/8$, as indicated by result (2.2),
and the convergence is by slow monotone decrease. 
% \fermat{finish this, Nils}
In fact, somewhat long calculations involving Taylor expansions
and moments of the $\chi^2_{n-1}$ yield  
$$\mise=n^{-1}\sigma^{-1}(2\sqrt{\pi})^{-1}
	\{7/8+(271/96)/n+O(1/n^2)\}. $$
% I found 7/8 + (271/96)/n + (14149/1024)/n^2 ... 

\section
{\sl 2.4. The UMVU density estimator.}
There is a unique uniformly minimum variance unbiased estimator 
of $f(x,\theta)=\sigma^{-1}\phi(\sigma^{-1}(x-\mu))$ under model conditions.
We include this in our list of competitors too:
$$\bar f(x)={k_n\over \hatt\sigma}\Bigl\{1-{n\over (n-1)^2}
	\Bigl({x-\hatt\mu\over \hatt\sigma}\Bigr)^2\Bigr\}^{n/2-2}
	\quad {\rm where\ }k_n=
{\Gamma(\half(n-1))\over 
	\sqrt{\pi}\Gamma(\half(n-2))}{n^{1/2}\over n-1}, \eqno(2.7)$$
in the region where $|x-\hatt\mu|/\hatt\sigma\le (n-1)/n^{1/2}$.
Outside this (random) region the density estimate is zero.
See Lehmann (1983, p.~87) for a sketch of a derivation
(where the same estimator is given in a slightly different form). 

A natural question is the relative comparison between the 
direct parametric estimator (1.1), which uses UMVU estimates 
for the parameters, and the UMVU density estimator $\bar f$. 
It is difficult to obtain a useful expression for the 
pointwise \mse{} of $\bar f$; 
stochastic simulation or two-dimensional numerical integration
seem to be called for. 
However, with some work the exact \mise{} can be obtained 
in a form suitable for numerical evaluation. 
With a standard normal $f$, start out 
writing the integrated squared error as 
$$\ise={k_n^2\over \hatt\sigma}{n-1\over n^{1/2}}
	\int_{-1}^1(1-v^2)^{n-4}\,\d v
	-2\int \bar f(x)f(x)\,\d x+R(f). $$
By the unbiasedness property of $\bar f$ this leads to 
$$\mise=\Bigl({n-1\over 2}\Bigr)^{1/2}
{\Gamma(\half(n-2))\over \Gamma(\half(n-1))}
	k_n^2{n-1\over n^{1/2}}
	{\Gamma(n-3)\sqrt{\pi}\over \Gamma(\half(2n-5))}
	-{1\over 2\sqrt{\pi}}\,. \eqno(2.8)$$
This can now be compared to the size of (2.6) for various $n$. 
It turns out that the UMVU density estimator actually
loses, for each $n$, to the simpler plug-in method. 
The ratio of (2.8) by (2.6) is always greater than 1,
but by an amount that quickly becomes small with $n$,
see the table in Section 5. 
% The ratio is 1.038 for $n=10$, 1.018 for $n=15$,
% 1.011 for $n=20$, and 1.003 for $n=50$, for example. 
That the ratio tends to 1 also reflects the fact that 
the (2.7) estimator is quite close to the plug-in (1.1)
estimator for moderate to large $n$. 

\bigskip
{\bf 3. MSE and MISE for the nonparametric estimator.}
This section works out the necessary exact expressions
for the bias and variance of the kernel estimator (1.2). 

\section
{\sl 3.1. Large-sample approximations.}
Again our focus is on the finite-sample case to follow afterwards,
but it is illuminating to compare those results to 
the traditional approximations suggested by well-known 
asymptotics arguments. 

The leading terms of the bias and variance are respectively 
$\half k_2h^2f''(x)$ and $(nh)^{-1}\allowbreak R(K)f(x)-n^{-1}f(x)^2$,
see for example Scott (1992, Ch.~6).
Here $k_2=\int u^2K(u)\,\d u$ while 
$R(q)$ is generic notation for $\int q(x)^2\,\d x$.
Hence the \mise{} can be approximated by 
$(nh)^{-1}R(K)+{1\over 4}k_2^2h^4R(f'')-n^{-1}R(f)$,
as $h$ goes to zero with $nh$ tending to infinity. 
The minimiser of this approximate \mise{} is 
$h_a=\{R(K)/k_2^2\}^{1/5}R(f'')^{-1/5}/n^{1/5}$, 
leading to a minimum approximate \mise{} 
of the form $d_n/n^{4/5}-R(f)/n$, where 
$d_n={5\over 4}\{R(K)k_2^{1/2}\}^{4/5}R(f'')^{1/5}$. 
\eject

These formulae can be suitably translated 
to the situation where $f$ is normal, for our two choices 
$K=\phi$ and $K=\Yep$. For the normal kernel case,
$$h_a\doteq 1.0592\,\sigma/n^{1/5} \quad {\rm and} \quad 
	\min_{h}\mise\doteq 0.3329\,\sigma^{-1}/n^{4/5}
	-0.2821\,\sigma^{-1}/n. $$
For the \Yepa{} kernel case, 
$$h_a\doteq 4.6898\,\sigma/n^{1/5} \quad {\rm and} \quad
	\min_{h}\mise\doteq 0.3198\,\sigma^{-1}/n^{4/5}
	-0.2821\,\sigma^{-1}/n. $$

\section
{\sl 3.2. Exact finite-sample \mse{} performance.} 
The mean and variance of the (1.2) estimator can 
generally be expressed as respectively 
$$e_h(x)=\E K_h(X_i-x)=\int K(u)f(x+hu)\,\d u \eqno(3.1)$$
and 
$$v_{n,h}(x)^2=n^{-1}\Var K_h(X_i-x)
	=(nh)^{-1}a_h(x)-n^{-1}e_h(x)^2, \eqno(3.2)$$
where $a_h(x)=\int K(u)^2f(x+hu)\,\d u$. 
Our task involves computing, displaying and comparing 
these function when $f$ is normal. It suffices to
take care of the standard normal case, giving 
certain $e_h^0(x)$ and $a_h^0(x)$ functions, say, 
in that the case of a normal $(\mu,\sigma^2)$ density 
then will lead to respectively 
$$e_h(x)=\sigma^{-1}e^0_{h/\sigma}(y)
	\quad {\rm and} \quad 
  a_h(x)=\sigma^{-2}a^0_{h/\sigma}(y), \eqno(3.3)$$ 
in terms of $y=(x-\mu)/\sigma$.

Consider first the case $K=\phi$. We there find 
$$\eqalign{
e_h^0(x)&=\phi\Bigl({x\over (1+h^2)^{1/2}}\Bigr){1\over (1+h^2)^{1/2}}, \cr
a_h^0(x)&={1\over 2\sqrt{\pi}}\phi\Bigl({x\over (1+h^2/2)^{1/2}}\Bigr)
	{1\over (1+h^2/2)^{1/2}}. \cr}$$
We may now compute $\mse(x)$ for any $x$, via formulae (3.1)--(3.3). 

More troublesome is finding the \mse{} for the \Yepa{} kernel
$\Yep$. It is convenient to define 
$$N_{j,h}(x)=h^{-1}\int_{x-h/2}^{x+h/2} v^j\phi(v)\,\d v
	\quad {\rm for\ }j=0,1,2,3,4. $$
Note that $N_{j,h}(x)=x^j\phi(x)+O(h^2)$ for small $h$. 
Values of these can be found exactly in terms of the 
$\phi$ and $\Phi$ functions, the density and cumulative 
of the standard normal, in that 
$\int_a^b v\phi(v)\,\d v=\phi(a)-\phi(b)$, 
$\int_a^b v^2\phi(v)\,\d v=a\phi(a)-b\phi(b)+\Phi(a,b)$, 
while 
$$\int_a^b v^j\phi(v)\,\d v=a^{j-1}\phi(a)-b^{j-1}\phi(b)
	+(j-1)\int_a^b v^{j-2}\phi(v)\,\d v
	\quad {\rm for\ }j\ge 3. \eqno(3.4)$$
These formulae are found by simple partial integration. Thus 
$$\eqalign{e^0_h(x)
&=(3/2)h^{-1}\int_{x\pm h/2}\Bigl\{1-4\Bigl({v-x\over h}\Bigr)^2\Bigr\}
	\phi(v)\,\d v \cr
&=(3/2)\bigl[N_{0,h}(x)-4h^{-2}\{N_{2,h}(x)
	-2xN_{1,h}(x)+x^2N_{0,h}(x)\}\bigr] \cr}$$
and similarly 
$$\eqalign{a^0_h(x)
&={9\over 4}h^{-1}\int_{x\pm h/2}\Bigl\{1-4\Bigl({v-x\over h}\Bigr)^2\Bigr\}^2
	\phi(v)\,\d v \cr
&={9\over 4}\bigl[N_{0,h}(x)
	-8h^{-2}\{N_{2,h}(x)-2xN_{1,h}(x)+x^2N_{0,h}(x)\} \cr
&+16h^{-4}\{N_{4,h}(x)-4xN_{3,h}(x)+6x^2N_{2,h}(x)
	-4x^3N_{1,h}(x)+x^4N_{0,h}(x)\}\bigr]. \cr}$$
Again, the exact bias, variance and \mse{} can now be computed for any $x$,
via (3.1)--(3.3). 

\section
{\sl 3.3. Exact finite-sample \mise{} calculations.} 
An exact formula for the \mise{} for the case of 
a normal kernel and a normal estimand is not difficult to obtain,
and was first exhibited by Fryer (1976). 
% see also Silverman (1986, Section 3.2), for example. 
It also follows below via other arguments. 
% [xx can be dropped: 
% Later far-reaching generalisations can be found 
% in Marron and Wand (1992) and Hjort and Glad (1996). xx] 
Presently we also need such a formula when $K=\Yep$, 
which is a more demanding technical exercise. 
In general it turns out to be helpful to work in terms of 
$g(y)=\int f(x)f(x+y)\,\d x$, the density of a difference $X_i-X_j$,
and with the corresponding $g_K(y)=\int K(u)K(u+y)\,\d u$. 
By suitably interchanging the order of integration in 
a couple of multiple integrals, using Fubini's theorem, one finds that 
$\int a_h(x)\,\d x=R(K)$, that 
$\int e_h(x)^2\,\d x=\int g_K(u)g(hu)\,\d u$, and that 
$\int e_h(x)f(x)\,\d x=\int K(u)g(hu)\,\d u$. Hence 
$$\mise(h)=(nh)^{-1}R(K)+(1-n^{-1})\int g_K(u)g(hu)\,\d u
	-2\int K(u)g(hu)\,\d u+g(0) \eqno(3.5)$$
is a generally valid formula for the exact \mise, with
$g(0)=R(f)$. 

For $K=\phi$, $g(u)$ becomes $\phi(u/\sqrt{2})/\sqrt{2}$, 
and the calculations become easy when $f$ is standard normal. 
They give 
$$\mise^0(h)=(2\sqrt{\pi})^{-1}
	\{(nh)^{-1}+(1-n^{-1})(1+h^2)^{-1/2}
	 	-2(1+\half h^2)^{-1/2}+1\}, \eqno(3.6)$$
where the `0' superscript indicates that the $f$ in question is $N(0,1)$.
For a general $N(\mu,\sigma^2)$ one finds 
$$\mise(h)=\sigma^{-1}\mise^0(\sigma^{-1}h). \eqno(3.7)$$
This can be minimised over $h$ for each sample size $n$
to find the optimal possible performance; see Section 4. 

Next turn to the $\Yep$ kernel, for which answers 
require more efforts. Rather than integrating 
the exact \mse{} expressions found above over $x$,
which is possible but rather involved, we prefer using the (3.5)
identity. It requires finding $g_K$ for the $\Yep$ kernel.
One sees that $g_K(u)$ is zero outside $[-1,1]$. 
For $0\le u\le 1$ somewhat long calculations give 
$$\eqalign{g_K(u)
&=\int K(x-\half u)K(x+\half u)\,\d x \cr
&=\int_{-(1-u)/2}^{(1-u)/2}
	\{1-4(x-\half u)^2\}\{1-4(x+\half u)^2\}\,\d x \cr 
&=(6/5)(1-5u^2+5u^3-u^5), \cr}$$
while $g_K(u)=(6/5)(1-5u^2-5u^3+u^5)$ for $-1\le u\le0$,
by symmetry. Let again $f$ be standard normal;
the (3.7) scale identity also holds for the $\Yep$ kernel. 
Further calculations give that $\int K(u)g(hu)\,\d u$ is equal to 
$$\eqalign{{\rm I}(h)
&=h^{-1}\int K(\sqrt{2}v/h)\phi(v)\,\d v \cr
&=h^{-1}3\int_0^{h/2\sqrt{2}}(1-4\cdot2v^2/h^2)\phi(v)\,\d v \cr
&=3h^{-1}\{(1-8/h^2)\Phi(0,h/2\sqrt{2})
	+2\sqrt{2}h^{-1}\phi(h/2\sqrt{2})\}. \cr}$$
Similarly, $\int g_K(u)g(hu)\,\d u$ is equal to 
$$\eqalign{{\rm II}(h)
&=h^{-1}\int g_K(\sqrt{2}v/h)\phi(v)\,\d v \cr
&=h^{-1}(12/5)\int_0^{h/\sqrt{2}}
	(1-5\cdot2v^2/h^2+5\cdot2\sqrt{2}v^3/h^3
	-4\sqrt{2}v^5/h^5)\phi(v)\,\d v \cr
&=(12/5)h^{-1}[(1-10/h^2)\Phi(0,h/\sqrt{2})
	+(20\sqrt{2}/h^3-32\sqrt{2}/h^5)\phi(0) \cr 
&\qquad\qquad	
	+(\sqrt{2}/h-12\sqrt{2}/h^3+32\sqrt{2}/h^5)\phi(h/\sqrt{2})]. \cr}$$
This suffices, via (3.5), to produce an exact \mise{} formula 
for the kernel estimator with the $\Yep$ kernel:
$$\mise^0(h)=(6/5)(nh)^{-1}
	+(1-n^{-1})\,{\rm II}(h)-2\,{\rm I}(h)+1/(2\sqrt{\pi}). \eqno(3.8)$$
 
\bigskip
{\bf 4. The real MISE when a random bandwidth is used.}
We have not yet been explicit about which $h$ value to use 
in the kernel estimator. The results above give 
relevant guidelines and lead to natural data-based proposals.
This also necessitates some extra analysis about 
the real \mise. 

\section
{\sl 4.1. Deciding on the bandwidth.}
Knowing that the optimal size of the $h$ is $O(n^{-1/5})$,
we consider the \mise{} curves via bandwidths 
of the form $h=b/n^{1/5}$, focusing on the $b$ parameter. 
For the normal kernel case (3.6) leads to 
$$\mise^0\Bigl({b\over n^{1/5}}\Bigr)
={1\over 2\sqrt{\pi}}
\Bigl[{1\over bn^{4/5}}+\Bigl(1-{1\over n}\Bigr)
	{1\over (1+b^2/n^{2/5})^{1/2}}
	-{2\over (1-\half b^2/n^{2/5})^{1/2}}+1\Bigr]. \eqno(4.1)$$
According to this and (3.7), 
the best theoretical bandwidth to use for estimation of 
a normal estimand is $h_{\rm mise}=b_n\sigma/n^{-1/5}$,
where $b_n$ is the minimiser of (4.1). 
The theory reviewed in Section 3.1 predicts 
that the minimum should occur around $b=1.059$ when $n$ is large. 
Some exact values are given in the table of Section 5.
One should actually stretch the bandwidth a little bit 
in comparison with the famous $1.059\,\sigma/n^{1/5}$ rule of thumb,
thereby earning a couple of percent in \mise. 
The actual rule we use in practice is $\tilda h=b_n\hatt\sigma/n^{1/5}$. 

Similarly, for the \Yepa{} kernel case, 
consider bandwidths of the form $h=c/n^{1/5}$, with 
$$\mise^0(c/n^{1/5})
=(6/5)c^{-1}n^{-4/5}+(1-n^{-1})\,{\rm II}(cn^{-1/5})-2\,{\rm I}(cn^{-1/5})
	+(2\sqrt{\pi})^{-1}, \eqno(4.2)$$
in which 
$${\rm I}\Bigl({c\over n^{1/5}}\Bigr)
=3{n^{1/5}\over c}\Bigl[\Bigl(1-{8n^{2/5}\over c^2}\Bigr)
	\Phi\Bigl(0,{c\over 2\sqrt{2}n^{1/5}}\Bigr)
	+2\sqrt{2}{n^{1/5}\over c}
	\phi\Bigl({c\over 2\sqrt{2}n^{1/5}}\Bigr)\Bigr], $$
$$\eqalign{{\rm II}\Bigl({c\over n^{1/5}}\Bigr)
&={12\over 5}{n^{1/5}\over c}\Bigl[\Bigl(1-{10n^{2/5}\over c^2}\Bigr)
	\Phi\Bigl(0,{c\over \sqrt{2}n^{1/5}}\Bigr)
+\Bigl({20\sqrt{2}n^{3/5}\over c^3}-{32\sqrt{2}n\over c^5}\Bigr)\phi(0) \cr
&\qquad\qquad 
+\Bigl({\sqrt{2}n^{1/5}\over c}-{12\sqrt{2}n^{3/5}\over c^3}
	+{32\sqrt{2}n\over c^5}\Bigr)
	\phi\Bigl({c\over \sqrt{2}n^{1/5}}\Bigr)\Bigr]. \cr}$$
The theoretically best bandwidth for a normal estimand 
is $h_{\rm mise}=c_n\sigma/n^{1/5}$, where $c_n$ minimises (4.2). 
The rule we use here is therefore $\tilda h=c_n\hatt\sigma/n^{1/5}$. 
By results of Section 3.1, $c_n$ should be close to 4.690 for large $n$. 
Exact values are provided in Section 5's table for 
some sample sizes. The convergence of $c_n$ towards its 
limiting value 4.690 is slow and from above. 

\section
{\sl 4.2. Determining the real MISE.} 
We used formulae (4.1) and (4.2) to decide on methods 
for choosing $h$, but as the $\tilda h$ used depend on data 
we cannot use the \mise{} formulae to evaluate the 
real \mise{} incurred. We should consider the 
expected value of the integrated squared error of the real estimator,
that is, 
$$\mise=\E\,\ise(\tilda f(\cdot,\tilda h))
	=\E\int\{\tilda f(x,\tilda h)-f(x)\}^2\,\d x. \eqno(4.3)$$
Jones (1991) calls this {\csc Eise} to make the distinction clear;
see his paper for a general discussion of the performance criteria used. 
If $\hatt\sigma$ is reasonably precise for $\sigma$, 
then $\tilda h$ is in the vicinity of $h_{\rm mise}$,
and the $\mise(h)$ curves are reasonably flat around their minimands. 
Hence the additional uncertainty in $\tilda h$ is bound to 
be a secondary effect, not likely to give a 
significant increase in the \mise{} level compared to 
$\mise(h_{\rm mise})$, at least when $n$ is moderate or big. 
However, precise numbers are called for in order to 
finish our comparison study for small $n$. 
Two methods are now described for obtaining these numbers.
We still assume that the true density is normal,
and may take it to be in standard form. 

The first method is exact, but requires numerical integration. 
Start out writing the \ise{} of a kernel estimator as 
$$\eqalign{
&{1\over n^2}\sum_{i,j}\int
	K_{\tilda h}(x-X_i)
	K_{\tilda h}(x-X_j)\,\d x 
	-2{1\over n}\sumin\int K_{\tilda h}(x-X_i)
		f(x)\,\d x+R(f) \cr
&={R(K)\over n\tilda h}+{1\over n^2}\sum_{i\not=j}
	{1\over \tilda h}g_K\Bigl({X_i-X_j\over \tilda h}\Bigr)
	-2{1\over n}\sumin\int K(u)f(X_i+\tilda hu)\,\d u+R(f), \cr}$$
which has mean value 
$$\mise=\E{R(K)\over n\tilda h}+\Bigl(1-{1\over n}\Bigr)
	\E{1\over \tilda h}g_K\Bigl({X_1-X_2\over \tilda h}\Bigr)
	-2\,\E\int K(u)f(X_1+\tilda hu)\,\d u+R(f). $$
This is valid if only the $\tilda h$ operation is symmetric in
the $n$ data points. Now focus on $\tilda h=a\hatt\sigma$ 
for a known value of $a$, like for the rules described above. 
Let $p_n(r)$ and $q_n(s)$ be the densities of respectively 
$R=(X_1-\hatt\mu)/\hatt\sigma$ and $S=(X_1-X_2)/\hatt\sigma$. 
These are ancillary statistics and as such are independent 
of $(\hatt\mu,\hatt\sigma)$, by the so-called Basu Lemma,
cf.~Lehmann (1983, Ch.~1.5). Hence 
$$\eqalign{\mise
&=R(K)(na)^{-1}\E\,\hatt\sigma^{-1}
	+(1-n^{-1})\E\,\hatt\sigma^{-1}\,\E a^{-1}g_K(a^{-1}S) \cr 
&\qquad\qquad 
	-2\,\E\int K(u)f(\hatt\mu+\hatt\sigma(R+au))\,\d u
		+R(f), \cr}\eqno(4.4)$$
and this can be evaluated exactly. 
The inverse mean of $\hatt\sigma$ was mentioned in Section 2.3 
% $\E\,\hatt\sigma^{-1}=\{\half(n-1)\}^{1/2}
% \Gamma(\half(n-2))/\Gamma(\half(n-1))$ 
and $\E a^{-1}g_K(a^{-1}S)=\int a^{-1}g_K(a^{-1}s)q_n(s)\,\d s$. 
Furthermore, for the $-2$ term, some calculations give 
$$\eqalign{
\E&\int\int K(u)f(\hatt\mu+\hatt\sigma(r+au))p_n(r)\,\d r\,\d u \cr 
&=\int\int K(u){1\over (2\pi)^{1/2}}\Bigl({n\over n+1}\Bigr)^{1/2}
	\Bigl\{1+{n\over n+1}{(r+au)^2\over n-1}\Bigr\}^{-\half(n-1)}
	p_n(r)\,\d r\,\d u. \cr}$$
Finally needed are the densities of $R$ and $S$, 
and with the required stamina one finds 
$p_n(r)=k_n\{1-nr^2/(n-1)^2\}^{\half(n-4)}$ for $|r|\le (n-1)/n^{1/2}$,
where $k_n$ is as in (2.7), and 
$q_n(s)=\lambda_n\{1-s^2/2(n-1)\}^{\half(n-4)}$ for $|s|\le \{2(n-1)\}^{1/2}$,
with $\lambda_n=k_n(n-1)^{1/2}/(2n)^{1/2}$. 

The above is sufficient to evaluate the real \mise{}
for the two kernel estimators, as described above, 
with $a=b_n/n^{1/5}$ for $K=\phi$ and $a=c_n/n^{1/5}$ for the \Yepa{} kernel. 

Our second method is simpler to execute but relies on
simulation and computer power, and is hence not quite as 
good as the exact solution. Write the \mise{} as 
$$\eqalign{\mise
% &=\E\int\{\tilda f(x,\tilda h)-f(x)\}^2\,\d x \cr 
&=\E\int f(x)\{\tilda f(x,\tilda h)/f(x)^{1/2}-f(x)^{1/2}\}^2\,\d x \cr 
&=\E_{n+1}\{\tilda f(X_{n+1},\tilda h)/f(X_{n+1})^{1/2}
		-f(X_{n+1})^{1/2}\}^2, \cr}$$
where the $\E_{n+1}$ indicates expectation with respect to 
$n+1$ variables $X_1,\ldots,X_{n+1}$ independently drawn from $f$;
of course $X_1,\ldots,X_n$ are present in the construction 
of $\tilda f(x,\tilda h)$. 
The great advantage is that a numerical value for $\mise$ 
can now be determined by simulation, 
without any integration, regardless of the complexities 
of the data-dependent rule used for $\tilda h$. 
It is also possible and sometimes advantageous, 
with respect to computer time and required precision, 
to evaluate the estimator constructed from
$X_1,\ldots,X_n$ in more than one new point, 
using say $\mise=\E_{n+m}Z_{n,m}$, where 
$$Z_{n,m}=m^{-1}\sum_{j=1}^m
	\{\tilda f(X_{n+j},\tilda h)/f(X_{n+j})^{1/2}
		-f(X_{n+j})^{1/2}\}^2. $$
The \mise{} is obtained by simulating a large enough number of 
such variables. The accuracy is essentially determined by
the size of $B$, the number of times the $(X_1,\ldots,X_n)$ set 
is sampled to form the density estimator, 
but is also helped by $m$, the number of sampled 
extra points at which to evaluate the constructed estimator. 
Some calculations show that the variance of the simulation based 
\mise{} figure has the approximate form 
$\{\tau_1^2/(Bm)+\tau_2^2/B\}/n^{8/5}$, where $\tau_1$ and $\tau_2$ 
are certain constants. The \mise{}s themselves go down
to zero with rate $n^{-4/5}$, and we need precision on the 
relative and not absolute scale, so the point is to get 
a low $\tau_1^2/(Bm)+\tau_2^2/B$. 

In our calculations both methods were employed. 
In the summary below only the exact numbers,
arrived at using numerical integration, are given.  
These have a clear edge in precision; 
even with a million simulations per sample size one cannot
be fully confident about the third decimal place 
when giving \mise{} numbers on appropriate relative scale. 

\bigskip
{\bf 5. Comparisons and illustrations.}
This section offers numerical results and two figures 
to illustrate our findings. 

The table provides full \mise-comparison information
about the two parametric and two nonparametric competitors
under discussion, for a selection of sample sizes. 
It is convenient to view the \mise{} obtained by the 
parametric method as a benchmark, then dividing 
other \mise{} numbers with this benchmark, 
to avoid staring at numbers that all go to zero. 
Such a ratio below one means outperforming the direct
parametric plug-in method, and vice versa. 
One immediate fact to note is that the parametric UMVU method 
always loses to the plug-in method, but with a relative 
amount that becomes smaller with $n$. The \mise{} for
the UMVU method for $n=3$ is actually infinite. 

From the table we see that the kernel method with a
normal kernel wins over the parametric method 
for $n\le 14$, if the optimal bandwidth is available,
and with the random bandwidth used in practice the 
verdict is $n\le 9$. Corresponding findings for 
the \Yepa{} kernel method are $n\le 15$ and $n\le 10$. 
The table in particular sheds light on the extra 
estimation noise for the kernel methods caused by 
having to estimate the best bandwidth from data. 

Also of general interest are the implicit finite-sample corrections
to the well-known `normal reference rules' for bandwidth
selection. Rather than using $1.059\allowbreak\,\hatt\sigma/n^{1/5}$
with a normal kernel one should stretch it further to the appropriate 
$b_n\,\hatt\sigma/n^{1/5}$; this will lead to a small
reduction in \mise. 

{{\medskip\narrower\baselineskip11pt
\parindent0pt 
{\csc Table.} \sl 
The \mise{} behaviour when estimating a normal density 
is summarised here, for two parametric and two nonparametric estimators. 
For each sample size $n$ we give first the `benchmark' \mise{}
for the parametric plug-in estimator (column 2). 
The \mise{} numbers for the other methods are then given 
as ratios, dividing with the benchmark \mise. 
Column 3 gives this ratio for the UMVU estimator. 
Columns 5 and 6 give the \mise-ratio for the kernel method 
with the normal kernel, with respectively best possible and 
estimated best possible bandwidth.
Similarly columns 8 and 9 give these ratios for the kernel
method with the \Yepa{} kernel, 
again with best possible and estimated best possible bandwidth. 
Columns 4 and 7 display the best constants $b_n$ and $c_n$
involved in the bandwidth selection procedures for the 
two kernel estimators. 
\smallskip}}
\overfullrule=0pt 
{{\smallskip\obeylines\parindent0pt\baselineskip10pt\sl 
~~~~~~~~~~~~~parametric:~~~~~~~~~~~~~~~normal kernel:~~~~~~~~~~~~~Yepanechnikov kernel:
~~~~~$n$~~~plug-in~~~~UMVU~~~~~~~~$b_n$~~~~~~~ratio$_1$~~~~ratio$_2$~~~~~~~$c_n$~~~~~~~ratio$_1$~~~~ratio$_2$
\smallskip\tt 
~~~3~~0.23230~~~~~$\infty${\hskip0.01cm}~~~~1.2871~~0.208~~0.699~~~5.2822~~0.209~~0.727
~~~4~~0.11829~~1.5095~~~1.2628~~0.350~~0.727~~~5.2177~~0.349~~0.747
~~~5~~0.07969~~1.2110~~~1.2458~~0.459~~0.773~~~5.1737~~0.455~~0.786
~~~6~~0.06016~~1.1223~~~1.2331~~0.548~~0.824~~~5.1411~~0.541~~0.830
~~~7~~0.04835~~1.0822~~~1.2230~~0.623~~0.874~~~5.1156~~0.614~~0.874
~~~8~~0.04042~~1.0602~~~1.2148~~0.689~~0.922~~~5.0949~~0.677~~0.918
~~~9~~0.03472~~1.0466~~~1.2080~~0.748~~0.967~~~5.0776~~0.733~~0.960
~~10~~0.03044~~1.0375~~~1.2021~~0.801~~1.010~~~5.0628~~0.784~~0.9997
~~11~~0.02710~~1.0312~~~1.1970~~0.849~~1.050~~~5.0500~~0.830~~1.037
~~12~~0.02441~~1.0264~~~1.1925~~0.894~~1.088~~~5.0388~~0.872~~1.072
~~13~~0.02222~~1.0229~~~1.1885~~0.935~~1.124~~~5.0288~~0.911~~1.106
~~14~~0.02038~~1.0201~~~1.1849~~0.973~~1.157~~~5.0198~~0.948~~1.137
~~15~~0.01883~~1.0178~~~1.1816~~1.009~~1.189~~~5.0117~~0.982~~1.167
~~16~~0.01749~~1.0160~~~1.1786~~1.043~~1.220~~~5.0043~~1.015~~1.195
~~17~~0.01633~~1.0145~~~1.1759~~1.075~~1.249~~~4.9975~~1.045~~1.223
~~18~~0.01532~~1.0132~~~1.1734~~1.106~~1.276~~~4.9913~~1.074~~1.248
~~19~~0.01443~~1.0121~~~1.1711~~1.135~~1.303~~~4.9855~~1.102~~1.273
~~20~~0.01363~~1.0112~~~1.1689~~1.163~~1.328~~~4.9801~~1.128~~1.297
~~50~~0.00513~~1.0032~~~1.1368~~1.694~~1.824~~~4.8996~~1.631~~1.764
~100~~0.00252~~1.0014~~~1.1190~~2.150~~2.259~~~4.8540~~2.064~~2.175
1000~~0.00025~~1.0001~~~1.0842~~4.163~~4.214~~~4.7617~~3.983~~4.032
\smallskip}}

\fermat{For simplicity and convenience 
the figures are shown on pages 17 and 18}While the table
gives numbers related to an overall comparison,
Figures 1 and 2 give supplementary information on 
the pointwise level. Figure 1 conveys the essentials of our story. 
It compares the behaviour of the parametric plug-in 
and the kernel method with the \Yepa{} kernel, 
for sample size 14, when estimating an unknown normal density. 
We see that the kernel method has a much bigger bias than
the parametric one (the latter is almost zero), but that 
it nevertheless wins overall via its lower variance, 
in the central $x$-area. 
The figure also shows that the \mise-based findings are no
artifacts of the particular overall criterion used; 
the kernel method wins in terms of non-controversial 
mean squared error in the essential central region $(-1.53,1.53)$,
except in the small interval $(-0.17,0.17)$, and then only barely. 

Figure 2 is concerned with another theme, 
that of comparing two kernel estimators. 
Both of them are popular in practice, 
and perhaps in particular the normal kernel one
because of ease of computations and guaranteed smoothness. 
The only real claim of the \Yepa{} kernel choice is 
asymptotic in nature; it minimises both pointwise \mse{}
and overall \mise{} when $n$ is large. 
It does have a slight advantage for most finite $n$ too,
as illustrated here, outperforming its smoother sister 
very modestly in terms of \mse{} in the central area. 

Our calculations and results also make clear one point
which perhaps has been unanswered in the literature:
there cannot be a universal finite-sample theorem about
the optimality of the \Yepa{} kernel. The normal kernel
does better for $n=2$ and $n=3$, in terms of best 
possible performance, and indeed for all sample sizes up to $n=7$
in terms of performance with estimated bandwidths. 

\bigskip
{\bf 6. Concluding comments.} 

\section
{\sl 6.1. Other parametric families.} 
We have demonstrated that the simple-minded nonparametric method
can win, against odds, over the normal density estimator. 
The same would typically happen for other low-dimensional 
parametric families, and particularly so if the model is complex,
because of the sampling variability in the parameter estimates. 
Any model containing the normal, with three or more parameters, 
would produce density estimators that a fortiori 
would have larger \mse's and \mise's than estimator (1.1)
under normal assumptions, that is, the nonparametric (1.2) 
would win for even more sample sizes. 

As an example, consider the estimator based on 
a natural three-parameter skew-extended normal family, 
$$\hatt f(x)=\hatt\gamma
	\Phi\Bigl({x-\hatt a\over \hatt b}\Bigr)^{\hatt\gamma-1}
	\phi\Bigl({x-\hatt a\over \hatt b}\Bigr){1\over \hatt b}, $$
with estimates coming from maximum likelihood.
The motivation for this model is that $T_\gamma((X-a)/b)$
is standard normal for appropriate $(a,b,\gamma)$,
where $T_\gamma(y)=\Phi^{-1}(\Phi(y)^\gamma)$ 
is a natural skewness-inducing family of transformations. 
When $\gamma=1$ things are normal. 
In general terms, via arguments given in Section 2.1, 
$$n\,\mise{}\arr\int f(x,\theta)^2u(x,\theta)^\tr J(\theta)^{-1}
	u(x,\theta)\,\d x={\rm Tr}\{J(\theta)^{-1}L(\theta)\} $$
as $n$ grows, where 
$L(\theta)=\int f(x,\theta)^2u(x,\theta)u(x,\theta)^\tr\,\d x$. 
For this particular family some calculations show that 
the \mise{} is approximately $0.342/(n\sigma)$ for normal estimands, 
which is about 1.386 times bigger than the value $0.247/(n\sigma)$ 
for the normal density estimator (1.1). 
Extrapolating from a fuller version of the table of 
Section 5 leads one to the educated guess that 
the kernel methods will outperform the parametric density
estimator for sample sizes at least up to 30. 
% \fermat{thimk!}[xx Hence ... fill in with sample size guess xx] 
% see com11r etcetera. 
% J=[1, 0, c; 0, 2, -d; c, -d, 1];
% L=[0.5, 0, e; 0, f, -g; e, -g, h]/(2*sqrt(pi)). 

\section
{\sl 6.2. Shrinking the parametric estimator.} 
Observing that high variance is the main reason
for losing the competition, one might consider shrinking 
the parametric estimator, that is, using $k\hatt f(x)$ instead
for suitable $k$. This somewhat artificial trick can indeed 
give smaller \mise. We illustrate this working with the shrunk 
UMVU estimator $k\bar f(x)$, since calculations are a bit easier then. 
Its \mise{} becomes $k^2\mise+(k-1)^2R(f)$, in terms of 
the $\mise$ for the $\bar f$ estimator, and this is minimised for 
$k=R(f)/\{\mise+R(f)\}$. This is independent of $(\mu,\sigma)$
and so can be used for normal estimands, 
resulting in the smaller 
$$\mise_k=\mise\,R(f)/\{R(f)+\mise\}.$$ 
As \mise{} goes down as $1/n$ the effect is small, 
only improving the parametric performance slightly. 
% \fermat{thimk!}[xx Should include specific answer here. xx]

\section
{\sl 6.3. Other bandwidth selectors.} 
More ambitious $h$ choice methods abound in the 
density estimation literature, and most of them
will be more statistically variable than the 
normality-based rule-of-thumb type $a_n\hatt\sigma/n^{1/5}$ 
that we have employed. Thus higher exact \mise's than those 
arrived at in Section 4 will result. 
It is perfectly fair to use these rules-of-thumb here, however, 
when one aims at good performance for all normal estimands.

\section
{\sl 6.4. Other situations where nonparametric wins against odds.} 
We have worked with density estimators to make our 
basic educational point, since the contrast between 
the two camps is so pedagogically conspicuous there. 
However, there are other and simpler situations
where the same phenomenon, nonparametrics outperforming
parametrics, can be observed. 

Let $X_1,\ldots,X_n$ be i.i.d.~from a lognormal $(a,b^2)$ distribution,
and suppose it is required to estimate the mean $\mu$.
Then the parametrically based proposal is 
$\hatt\mu_{\rm pm}=\exp(\hatt a+\half\hatt b^2)$,
in that $\mu=\exp(a+\half b^2)$ under model conditions.
Here $\hatt a$ and $\hatt b$ are the natural estimators,
say $\hatt a=\bar A$ and $\hatt b^2=(n-1)^{-1}\sumin (A_i-\bar A)^2$,
in terms of $A_i=\log X_i$, and these are normal $(a,b^2)$.
But a competitor is the nonparametric $\hatt\mu_{\rm nonpm}=\bar X$.

The well-known asymptotic optimality result for maximum likelihood
estimation entails in particular that the parametric must be 
better than $\bar X$ for all large $n$,
and this is nicely illustrated here, in that calculations show 
$${\Var\,\bar X\over \Var\,\hatt\mu_{\rm pm}}
	\arr{\exp(b^2)-1\over b^2+\half b^4}
	={b^2+\half b^4+{1\over 6}b^6+\cdots
	  \over b^2+\half b^4}>1. $$
But for moderate to small $n$ the picture turns out 
to be different. We have 
$\mse(\bar X)=n^{-1}\exp(2a+b^2)\{\exp(b^2)-1\}$, 
while calculations give 
$$\mse(\hatt\mu_{\rm pm})=\exp(2a+b^2/n)\Bigl[\exp(b^2/n)
	\Bigl(1-{2b^2\over n-1}\Bigr)^{-(n-1)/2}
	-\Bigl(1-{b^2\over n-1}\Bigr)^{-(n-1)}\Bigr]. $$
Here one quickly notices that $b<\{\half(n-1)\}^{1/2}$ 
is required for the variance and \mse{} to exist at all 
for the parametric estimator, so for small $n$ the nonparametric
estimator is bound to win. This continues to be true 
also beyond the $n-1>2b^2$ requirement, as some analysis reveals. 
Thus for each $b$ there is a $n_0(b)$ such that the nonparametric
estimator wins over the parametric one, for each $n\le n_0(b)$,
while the parametric indeed eventually wins, namely for each $n>n_0(b)$.  
Values of $n_0(b)$ for $b=0.2,0.4,0.6,0.8,1.0,1.2$ are for example
found to be respectively 312, 87, 45, 31, 25, 22. 

\section
{\sl 6.5. Biasing can be worthwhile.} 
One may view our examples as part of 
a more general theme, which is that procedures 
that are biased but perhaps less variable 
sometimes are fruitful and strong competitors to 
the approximately unbiased ones. 
Cases in point include Stein estimation and 
more general shrinkage type methods. 
See also for example Hjort (1994, 1997) for analysis of procedures 
based on narrower versus wider parametric models.

\section
{\sl 6.6. Semiparametric attempts.} 
Our results might serve as a warning against putting too
much trust in intuition and large-sample results 
when judging the relative strengths of parametric and nonparametric methods. 
The findings may also have consequences for 
some of several recent semiparametric constructions 
that in various ways try to balance between the 
parametric and nonparametric frameworks. 
One relevant example is the mixture type estimator
$f^*(x)=\hatt p\,f(x,\hatt\theta)+(1-\hatt p)\,\tilda f(x)$,
worked with by Schuster and Yakowitz (1985),
Olkin and Spiegelman (1987) and others; 
more general methods with similar characteristics 
have been developed from a Bayesian perspective 
in Hjort (1996) and Byholt (1996). The crux is to have 
a good procedure for choosing $\hatt p$,
potentially also depending on $x$, from data. 
Intuition says perhaps that $\hatt p$ should be one or close to one
if the parametric model looks good, and indeed such methods have 
been proposed based on goodness of fit criteria. 
But this issue is more involved in view of the findings of 
this paper; even if the model is perfect we should not 
have $\hatt p$ close to one, for small sample sizes. 

\bigskip
\centerline{\bf References}

\def\ref#1{{\noindent\hangafter=1\hangindent=20pt
  #1\smallskip}}          
% my usual preference 
\parindent0pt
\baselineskip11pt
\parskip3pt 
\medskip 

\ref{%
Byholt, M. (1996).
Nonparametric density estimation with a local likelihood.
Cand.\allowbreak~scient.-thesis, Department of Mathematics,
University of Oslo.} 

% \ref{%
% Fan, J.~and Gijbels, I. (1996).
% {\sl Local Polynomial Modelling and its Applications.}
% Chapman and Hall, London.} 

% \ref{%
% Fenstad, G.U.~and Hjort, N.L. (1996).
% Statistical Reseach Report, Department of Mathematics,
% University of Oslo.}

\ref{%
Fryer, M.J. (1976).
Some errors associated with the non-parametric estimation 
of density functions. 
{\sl Journal of the Institute of Mathematics 
and its Applications}~{\bf 18}, 371--380.}

\ref{%
Hjort, N.L. (1994).
The exact amount of t-ness that the normal model can tolerate. 
{\sl Journal of the American Statistical Association} 
{\bf 89}, 665--675.}

\ref{%
Hjort, N.L. (1996). 
Bayesian approaches to semiparametric density estimation
(with discussion contributions).
In {\sl Bayes\-ian Statistics 5}, 
proceedings of the Fifth International Val\`encia Meeting 
on Bayesian Statistics (eds.~J.~Berger, 
J.~Bernardo, A.P.~Dawid, A.F.M.~Smith),
223--253. Oxford University Press.}

\ref{%
Hjort, N.L. (1997).
Estimation in moderately misspecified models. 
{\sl Journal of the American Statistical Association},
to appear.} 

% \ref{%
% Hjort, N.L.~and Glad, I.K. (1995).
% Nonparametric density estimation with a parametric start.
% {\sl Annals of Statistics} {\bf 23}, 882--904.}

% \ref{%
% Hjort, N.L.~and Glad, I.K. (1996).
% On the exact performance of a multiplicative 
% semiparametric density estimator. 
% Submitted for publication.} 

% \ref{%
% Hjort, N.L.~and Jones, M.C. (1996a).
% Locally parametric nonparametric density estimation.
% {\sl Annals of Statistics} {\bf 24}, to appear.} 

\ref{%
Jones, M.C. (1991).
The roles of ISE and MISE in density estimation.
{\sl Statistics and Probability Letters} {\bf 12}, 51--56.} 

\ref{%
Lehmann, E.L. (1983).
{\sl Theory of Point Estimation.}
Wiley, New York.}

% \ref{%
% Marron, J.S.~and Wand, M.P. (1992). 
% Exact mean integrated squared error. 
% {\sl Annals of Statistics} {\bf 20}, 712--736.}

\ref{%
Olkin, I.~and Spiegelman, C.H. (1987).
A semiparametric approach to density estimation.
{\sl Journal of the American Statistical Association} {\bf 82},
858--865.}

\ref{%
Schuster, E.~and Yakowitz, S. (1985).
Parametric/nonparametric mixture density estimation with application
to flood-frequency analysis. 
{\sl Water Resources Bulletin} {\bf 21}, 797--804.}

\ref{%
Scott, D.W. (1992).
{\sl Multivariate Density Estimation:
Theory, Practice, and Visualization.}
Wiley, New York.}

% \ref{%
% Silverman, B.W. (1986).
% {\sl Density Estimation for Statistics and Data Analysis.}
% Chapman and Hall, London.}

% \ref{%
% Simonoff, J. (1996).
% {\sl Smoothing in Statistics.}
% Springer-Verlag, New York.}

\ref{%
Wand, M.P.~and Jones, M.C. (1995).
{\sl Kernel Smoothing.}
Chapman \& Hall, London.}

\bigskip
\bigskip
\parindent20pt
\baselineskip13pt
\parskip0pt 

\vfill\eject

% \smallskip
% \centerline{\psfig{file=morten_nils_1.ps,height=3.5in,width=5.5in,angle=360}}
% \vskip0.1truecm

\centerline{\includegraphics[scale=0.66]{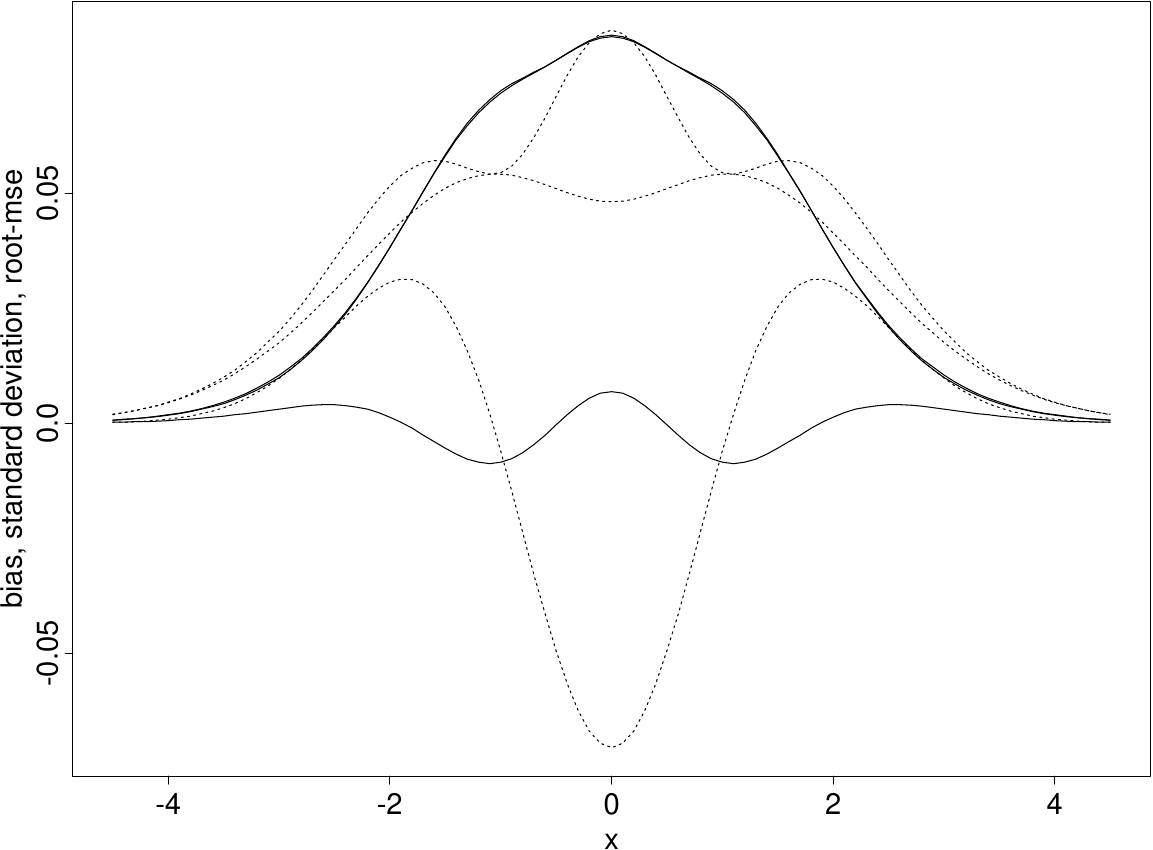}}

\bigskip\noindent
{\sl\narrower\noindent Figure 1: 
This figure illustrates behaviour of two density estimators 
when the estimand is standard normal and the sample size is 14. 
The three dotted lines display bias, standard deviation and 
root mean squared error for the kernel method (1.2) 
with the Yepanechnikov kernel. 
The three full lines similarly show bias, standard deviation 
and root mean squared error for the parametric estimator (1.1). 
The bias is quite small for the parametric method, 
which is why its root mean squared error 
is above its standard deviation with an almost invisible amount.  
  \smallskip}

\vfill\eject 

% \medskip\noindent
% \centerline{\psfig{file=morten_nils_2.ps,height=3.5in,width=5.5in,angle=0}}
% \vskip0.1truecm

\centerline{\includegraphics[scale=0.66]{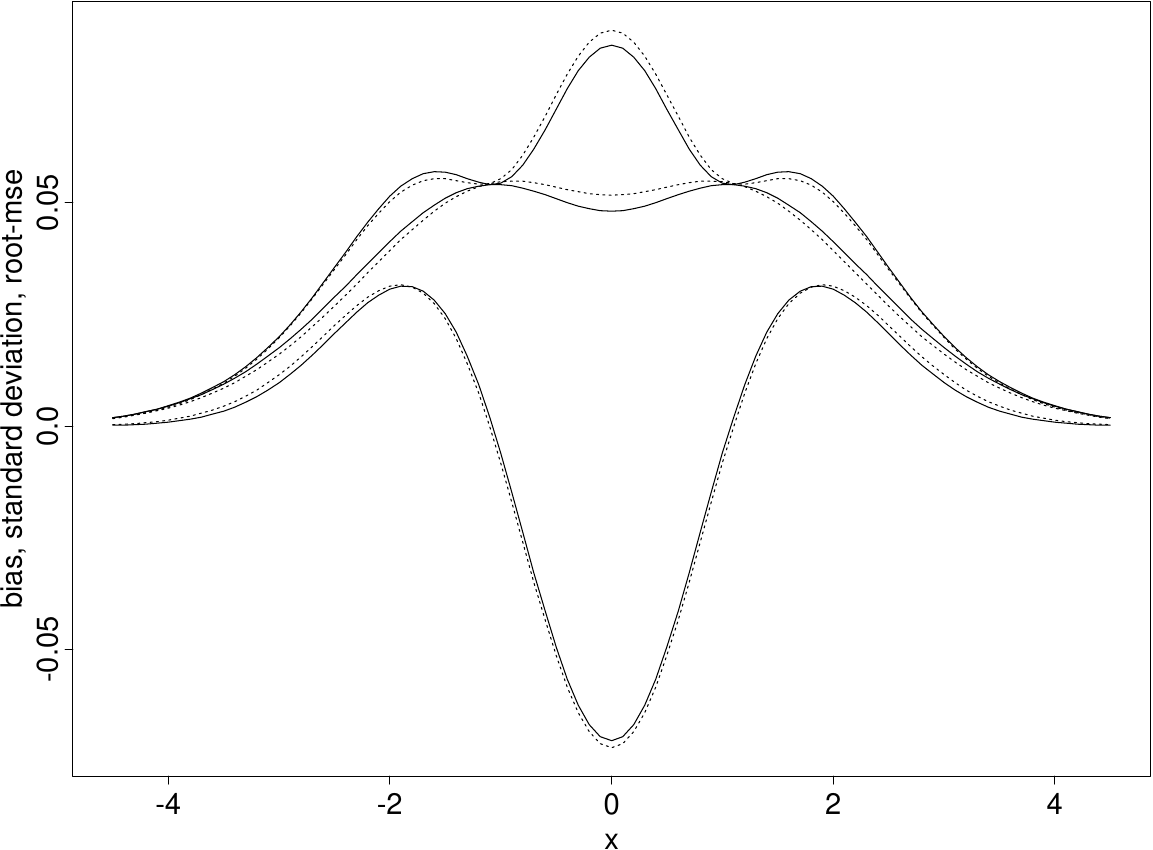}}

\bigskip\noindent
{\sl\narrower\noindent Figure 2:
This figure shows the relatively small difference in behaviour
between two competing kernel estimators, those using respectively
the normal kernel and the Yepanechnikov kernel.  
Displayed are bias, standard deviation and root mean squared error
for both methods, full lines for the Yepanechnikov kernel method 
and dotted lines for the normal kernel method. 
The estimand is standard normal, the sample size is 14,
and the best bandwidths are used. 
The Yepanechnikov kernel wins with about 
2.5 percent in terms of the overall criterion {\csc mise}. 
  \smallskip}

\vfill\eject 

\bye